\documentclass{ifacconf}

\usepackage{graphicx}      
\usepackage{amsmath,amssymb}
\usepackage[numbers]{natbib}        
\begin{document}
\begin{frontmatter}

\title{Hidden Markov chains and fields with observations in Riemannian manifolds} 


\author[First]{Salem Said} 
\author[Second]{Nicolas Le Bihan} 
\author[Third]{Jonathan H. Manton}

\address[First]{CNRS, University of Bordeaux (salem.said@u-bordeaux.fr)}
\address[Second]{CNRS, Gipsa-lab (nicolas.le-bihan@gipsa-lab.grenoble-inp.fr)}
\address[Third]{The University of Melbourne (jmanton@unimelb.edu.au)}

\begin{abstract}                
Hidden Markov chain, or Markov field, models, with observations in a Euclidean space, play a major role across signal and image processing. The present work provides a statistical framework which can be used to extend these models, along with related, popular algorithms (such as the Baum-Welch algorithm), to the case where the observations lie in a Riemannian manifold. It is motivated by the potential use of hidden Markov chains and fields, with observations in Riemannian manifolds, as models for complex signals and images.
\end{abstract}

\end{frontmatter}

\section{Introduction}
\vspace{-0.28cm}
 
The present work is concerned with hidden Markov chain and Markov field models, with observations in Riemannian manifolds.

It introduces a general formulation of hidden Markov chain models, whose observations lie in a Riemannian manifold, and derives an expectation-maximisation algorithm, to estimate the parameters of these models --- Sections \ref{sec:model} and \ref{sec:algorithm}, respectively.

It also describes a general hidden Markov field model, with observations in a Riemannian manifold, and discusses briefly the estimation of this model, using the expectation-maximisation approach --- Section \ref{sec:fields}. 

Hidden Markov chain and field models have been highly influential in signal and image analysis\!\cite{hmmbook}\!\!\cite{moulines}. For example,\linebreak the Baum-Welch algorithm, for estimating hidden Markov chain models, is a staple of applications such as speech recognition and protein sequencing\!\cite{rabiner}\!\!\!\cite{durbin}. Also, hidden Markov fields are central in several approaches to image restoration and segmentation\!\cite{stanli}. 

While quite extensive, existing literature is almost entirely focused on hidden Markov chain and field models with\hfill\linebreak observations in a Euclidean space. The present work aims to provide a statistical framework, which may be used to extend current algorithms, dealing with observations in a Euclidean space, to observations in more general Riemannian manifolds. 

For example, the expectation-maximisation algorithm,\hfill\linebreak derived in Section \ref{sec:algorithm}, extends the popular Baum-Welch algorithm to hidden Markov chains with observations in a Riemannian manifold.

The motivation for considering hidden Markov chains and fields with observations in Riemannian manifolds\hfill\linebreak comes from the important role which data in Riemannian manifolds can play in the study of complex signals and\hfill\linebreak images (or even three-dimensional fields, such as turbulent flows\!\cite{color}). 

To understand this, consider the case of a non-stationary, multivariate signal, observed over a long time interval. This is typical of a phased-array radar, or brain-computer interface signal\!\cite{jeuris}\!\!\cite{bci}. 

While this signal cannot be entirely described by a\hfill\linebreak single tractable model, it is still possible to ``zoom in" on short time windows, which admit of a complete, tractable description, very often in terms of a covariance matrix (for array radar signals, this has an additional block-Toeplitz structure\!\cite{jeuris}).  

Now, the original complex signal appears as a time series of covariance matrices, or other descriptors, one for each short time window, and it seems possible to extract useful information, by studying the statistical distribution of these descriptors. 

This idea can also be applied to a complex image, or three-dimensional field, by zooming in on small connected regions\!\cite{image-tuzel}\!\!\cite{color}. Recently, it has generated significant interest in the problem of estimating an unknown distribution,\hfill\linebreak underlying a population of descriptors (\textit{e.g.}, covariance matrices, principal subspaces, histogrammes)\!\cite{said2}\!\!\cite{chevallier}\!\!\cite{banerjee}. 

Mathematically, this problem was framed as the problem of estimating an unknown probability distribution on a Riemannian manifold. Indeed, in several applications, this Riemannian approach has become standard, due to its superior performance (\textit{e.g.}, brain-computer interface\!\cite{bci}).  

Existing work on this problem has always started from the assumption that descriptors, obtained from different time windows of a signal (or regions of an image), are sampled independently, all from the same underlying distribution. \hfill\linebreak
This assumption obscures the dynamics of the signal (or the spatial dependence structure within the image), potentially leading to the loss of crucial information. 

The statistical framework introduced in the present work captures this information, in the form of a hidden Markov dependence structure, which controls the distribution of the descriptors.

\section{Hidden Markov chains} \label{sec:model}
Hidden Markov chains are popular models for many kinds of signals (neuro-biological, speech, \textit{etc.}\!\cite{hmmbook}\!\!\cite{moulines}). Here, these models are generalised to manifold-valued signals.

In a hidden Markov chain model, one is interested in understanding some hidden, time-varying, finitely-valued, process, say $(q_t\,;t=1,2,\ldots)$ which takes values in a finite set $S$. When $q_t = \alpha$ for some $\alpha \in S$, one says that the process is in state $\alpha$ at time $t$.

Moreover, it is assumed that $(q_t)$ is a time-stationary Markov chain, so there exists a so-called transition matrix, $(P_{\alpha\beta}\,; \alpha,\beta \in S)$, which gives the conditional probabilities
\begin{equation} \label{eq:mchain}
  \mathbb{P}(q_{t+1} = \beta|q_t = \alpha) \,=\, P_{\alpha\beta} 
\end{equation}
In particular, if $\pi_t(\alpha) = \mathbb{P}(q_t = \alpha)$ is the distribution of $q_t$, then\!\cite{norris}
\begin{equation} \label{eq:transition}
\pi_{t+1}(\beta) \,=\, \sum_{\alpha \in S}\, \pi_t(\alpha)P_{\alpha\beta} 
\end{equation}
describes the transition from time $t$ to time $t+1$.

The states $(q_t)$ are hidden, so they can only be observed through random outputs $(y_t)$ which have their values in a Riemannian manifold $M$. 

Moreover, these outputs are generated independently from each other, and each $y_t$ depends only on $q_t$\!\cite{hmmbook}\!\!\cite{moulines}.

It is assumed $M$ is a homogeneous Riemannian manifold, and $y_t$ is distributed according to a location-scale model, 
\begin{equation}\label{eq:locationscale}
  p(y_t|q_t = \alpha) \,=\, f(y_t|\bar{y}_{\alpha\,},\sigma_\alpha)
\end{equation}
where $p(y_t|q_t = \alpha)$ denotes the conditional density, with respect to the Riemannian volume measure of $M$. 

The location parameters $\bar{y}_\alpha \in M$, and scale parameters $\sigma_\alpha \geq 0$ are unknown, but the function $f(y_t|\bar{y}_{\alpha\,},\sigma_\alpha)$ is known, and of the form\!\cite{warped}
\begin{equation} \label{eq:f}
f(y|\bar{y},\sigma) \,=\,\exp\left[\, \eta(\sigma)D(y\,,\bar{y}) - \psi(\eta(\sigma))\right]
\end{equation}
where $D:M\times M \rightarrow \mathbb{R}$, and where $\eta(\sigma)$ is the so-called natural parameter, while $\psi(\eta)$ is a strictly convex function. 

Examples of location-scale models of the form (\ref{eq:f}) include the von Mises-Fisher model, and the Riemannian Gaussian model, which are briefly recalled, in Remark 1, at the end of the present section.

From the definition of $\pi_t$ and from (\ref{eq:locationscale}), it follows that the probability density of $y_t$ is a time-varying mixture density
\begin{equation} \label{eq:density}
  p(y_t) \,=\, \sum_{\alpha \in S}\,\pi_t(\alpha)\,f(y_t|\bar{y}_{\alpha\,},\sigma_\alpha)
\end{equation}
Having access only to the observations $(y_t\,;t=1,2,\ldots)$, one hopes to recover as much information about the Markov chain $(q_t)$ as possible\,: its transition matrix $(P_{\alpha\beta})$ and the parameters $(\bar{y}_\alpha\,,\sigma_\alpha)$ which govern its outputs $y_t$.

The problem of estimating $(P_{\alpha\beta})$ and $(\bar{y}_\alpha\,,\sigma_\alpha)$ is addressed in Section \ref{sec:algorithm}. A further, equally interesting problem, will not be considered, for lack of space. This is the problem of estimating the invariant distribution of the chain $(q_t)$\!\cite{norris}\,: \\
if the transition matrix $(P_{\alpha\beta})$ is irreducible and aperiodic, then $\pi_t(\alpha) \rightarrow \pi(\alpha)$, the invariant distribution, as $t \rightarrow \infty$. This means that, in time, the observation $y_t$ will tend to sample from a mixture distribution of the same form (\ref{eq:density}), but with the invariant weights $\pi(\alpha)$ instead of  $\pi_t(\alpha)$.

$\diamond$ \textit{\textbf{Remark 1}\,:} location-scale models of the form (\ref{eq:f}) were introduced in\!\cite{warped}. For these models, $M$ is any homogeneous Riemannian manifold\!\cite{besse}\,: 
a Riemannian manifold with a group of isometries $G$ which acts transitively.

Isometry means each $g \in G$ is a mapping $g:M\rightarrow M$ which preserves the Riemannian metric of $M$. Moreover, transitive action means that for each $y,z \in M$ there exists $g \in G$ such that $g\cdot y = z$ (here, $g\cdot y = g(y))$.

The function $D : M \times M \rightarrow \mathbb{R}$, appearing in (\ref{eq:f}) is required to satisfy the group-invariance property
\begin{equation}
 D(y\,,\bar{y}) \,=\, D(g\cdot y\,,g\cdot \bar{y}) 
\end{equation}
for any $g \in G$ and $y,\bar{y} \in M$.

This property guarantees that $f(y|\bar{y},\sigma)$ is indeed a true probability density. In particular, the function $\psi(\eta)$ is strictly convex, because it is the cumulant generating function of the statistic $\Delta = D(y\,,\bar{y})$,
$$
\mathbb{E}_{(\bar{y},\sigma)}\exp(\eta\, \Delta)
\,=\,
e^{\psi(\eta)}
$$
where the expectation is with respect to the density $f(y|\bar{y},\sigma)$ -- this identity holds for any real $\eta$ as long as the expectation is finite. \hfill$\blacksquare$

Examples of location-scale models of the form (\ref{eq:f}) include the von Mises-Fisher model, and the Riemannian Gaussian model. \\[0.1cm]
$\diamond$ \textit{\textbf{von Mises-Fisher model}}\!\cite{mardia}\,: for this model, 
\begin{equation} \label{eq:vmf1}
M = S^{d-1} \hspace{0.2cm}\text{and}\hspace{0.2cm} G = O(d)
\end{equation}
where $S^{d-1} \subset \mathbb{R}^d$ is the $(d-1)$-dimensional unit sphere, and $O(d)$ is the group of orthogonal transformation of $\mathbb{R}^d$.

The function $D:M \times M \rightarrow \mathbb{R}$ is given by
\begin{equation} \label{eq:vmf2}
D(y\,,\bar{y}) \,=\, \langle y,\bar{y}\rangle_{\scriptscriptstyle \mathbb{R}^d}
\end{equation}
where $\langle \cdot,\cdot\rangle_{\scriptscriptstyle \mathbb{R}^d}$ is the Euclidean scalar product in $\mathbb{R}^d$, and $\psi(\eta)$ has the following expression (where $\nu = d/2$)
\begin{equation} \label{eq:vmf3}
  e^{\psi(\eta)} \,=\, (2\pi)^\nu\,\eta^{1-\nu}I_{\nu-1}(\eta)
\end{equation}
where $I_{\nu-1}$ is the modified Bessel function of order $\nu - 1$. Finally, the von Mises-Fisher model is obtained by setting $\eta(\sigma) = \sigma$ in (\ref{eq:f}) --- then, $\eta$ is called the concentration parameter. \hfill$\blacksquare$\\[0.1cm]
$\diamond$ \textit{\textbf{Riemannian Gaussian model}}\!\cite{said1}\!\!\cite{said2}\,: for this model
\begin{equation}\label{eq:rg1}
 M = \mathcal{P}_{\scriptscriptstyle d} \hspace{0.2cm}\text{and}\hspace{0.2cm} G = GL(d)
\end{equation}
where $\mathcal{P}_{\scriptscriptstyle d}$ is the space of $d \times d$ symmetric positive-definite matrices, and $GL(d)$ the group of invertible $d \times d$ matrices.

Here, $G$ acts on $M$ by $g\cdot y = gyg^\dagger$ where $^\dagger$ denotes transpose. Moreover, this action preserves the Riemannian distance, 
$$
d^2(y,z) \,=\, \mathrm{tr}\left[\left(\log(y^{-1}z)\right)^2\right]
$$
where $\mathrm{tr}$ denotes the trace, and $\log$ the symmetric matrix logarithm. The function $D:M \times M \rightarrow \mathbb{R}$ is then given by
\begin{equation} \label{eq:rg2}
D(y\,,\bar{y}) \,=\, d^2(y,\bar{y})
\end{equation}
and an expression of $\psi(\eta)$, in terms of a multivariate integral, is given in\!\cite{said1}. 

The model is obtained by setting $\eta(\sigma) = -\frac{1}{2\sigma^2}$ in (\ref{eq:f}). This gives the Riemannian Gaussian density
\begin{equation} \label{eq:rg3}
 f(y|\bar{y},\sigma) \,=\, Z^{-1}(\sigma)\exp\left[-\frac{d^2(y,\bar{y})}{2\sigma^2}\right]
\end{equation}
where $Z(\sigma) = e^{\psi(\eta(\sigma))}$. \hfill$\blacksquare$
\vfill
\pagebreak
\section{The EM algorithm} \label{sec:algorithm}
Here, an EM (expectation-maximisation) algorithm will be introduced, which addresses the estimation problem defined in Section \ref{sec:model}. 

To see why the expectation-maximisation approach is a natural choice, assume it were possible to access the so-called complete data $x_t =(y_{t\,},q_t)$ at times $t = 1,\ldots,T$. The resulting log-likelihood function would be,
\begin{equation} 
\tiny{\ell(x|\theta) \,=\, \log\left[ \pi_{\scriptscriptstyle 1}(q_{\scriptscriptstyle1})\prod^{T-1}_{t=1}\,P_{q_tq_{t+1}}\right]+\log\left[\prod^T_{t=1} f(y_t|\bar{y}_{q_t\,},\sigma_{q_t})\right]}\\[0.2cm]
\end{equation}
where $\theta = (P_{\alpha\beta\,},\bar{y}_\alpha\,,\sigma_\alpha)$. However, this could also be written\footnote{the term depending on $\pi_{\scriptscriptstyle 1}$ is discarded, since it is statistically non-significant when $T$ is large.},
\begin{equation} \label{eq:lc}
\tiny{\ell(x|\theta) = \sum_{\alpha,\beta \in S}N_{\alpha\beta}(q)\log\left(P_{\alpha\beta}\right) \,+\, \sum_{\alpha \in S}\sum^T_{t=1} \delta_\alpha(q_t)\log\,f(y_t|\bar{y}_{\alpha\,},\sigma_{\alpha})}
\end{equation}
where $N_{\alpha\beta}(q)$ is the number of times that $(q_{t\,},q_{t+1}) =  (\alpha,\beta)$, and where $\delta_\alpha(q_t) = 1$ if $q_t = \alpha$ and $=0$ otherwise.

Estimating $\theta$ from the complete data, by maximising the log-likelihood function (\ref{eq:lc}), is relatively straightforward, as explained in Remark 3 below. 

Unfortunately, since the Markov chain $(q_t)$ is hidden, one has access only to the observations $(y_t\,;t=1,\ldots,T)$. 

The log-likelihood function for these observations is much more complicated than (\ref{eq:lc}),
\begin{equation} \label{eq:lhmm}
 \tiny{
 \ell(y|\theta) = \log\left[
 \sum_{\scriptscriptstyle q_1 \in S}
\ldots
\sum_{\scriptscriptstyle q_T \in S} \pi_{\scriptscriptstyle 1}(q_{\scriptscriptstyle1})\prod^{T-1}_{t=1}P_{q_tq_{t+1}}\prod^T_{t=1} f(y_t|\bar{y}_{q_t\,},\sigma_{q_t})\right]
} \\[0.2cm]
\end{equation}
Instead of directly maximising (\ref{eq:lhmm}), the expectation-maximisation approach reduces the problem of estimating
$\theta$ from the observations $y_t$ to a sequence of steps which only involve maximising a function of the form (\ref{eq:lc}). 

Specifically, the $k$-th iteration of the EM algorithm computes estimates 
$$
\theta^k = (\hat{P}^k_{\alpha\beta\,},\bar{y}^k_\alpha\,,\sigma^k_\alpha)
$$
which are guaranteed to increase the log-likelihood (\ref{eq:lhmm}), in the sense that $\ell(y|\theta^k) \geq \ell(y|\theta^{k-1})$. Each iteration is made up of two steps\!\cite{dempster}\!\!\cite{moulines} \\[0.1cm]
$\bullet$ \textbf{E step}\,: compute the function
\begin{equation} \label{eq:estep}
  Q(\theta|\theta^{k-1}) \,=\, \mathbb{E}_{\theta^{k-1}}\left[\ell(x|\theta)\middle| y\right]
\end{equation}
where $\mathbb{E}_{\theta^{k-1}}$ means the expectation is computed under the current estimate $\theta^{k-1\,}$. This is called the expectation step. \\[0.1cm]
$\bullet$ \textbf{M step}\,: maximise $Q(\theta|\theta^{k-1})$ to obtain the new estimate
\begin{equation} \label{eq:mstep}
 \theta^{k} \,=\, \mathrm{argmax}_\theta\; Q(\theta|\theta^{k-1})
\end{equation}
In principle, the algorithm cycles through these two steps, until the estimates $\theta^k$ converge. In practice, this means the algorithm is designed to stop when the difference between $\theta^k$ and $\theta^{k-1}$ falls short of some pre-assigned accuracy.

There is no theoretical guarantee that the $\theta^k$ will converge to the global maximum of the log-likelihood (\ref{eq:lhmm}). In fact, convergence typically takes place only to a local maximum, depending on initialisation. 

\vfill
\pagebreak

However, the EM algorithm remains quite advantageous, from a computational point of view, since it boils down to relatively simple operations, with little computational complexity, in comparison with alternative methods.

Indeed, replacing (\ref{eq:lc}) into (\ref{eq:estep}), yields $Q(\theta|\theta^{k-1}) = Q(\theta)$, 
\begin{equation} \label{eq:qnew}
\tiny{Q(\theta) = \sum_{\alpha,\beta \in S}\nu_{\alpha\beta}(y)\log\left(P_{\alpha\beta}\right) \,+\, \sum_{\alpha \in S}\sum^T_{t=1} \omega^t_\alpha(y)\log\,f(y_t|\bar{y}_{\alpha\,},\sigma_{\alpha})} \\[0.15cm]
\end{equation}
with $\nu_{\alpha\beta}(y)$ and $\omega^t_\alpha(y)$ computed from the observations $y_t\,$,\hfill\linebreak
\begin{subequations} \label{eq:estepprac}
\begin{equation} \label{eq:nu}
\nu_{\alpha\beta}(y) \,=\, \sum^{T-1}_{t=1}\,\mathbb{P}_{\theta^{k-1}}(q_t = \alpha\,,q_{t+1}=\beta|y)
\end{equation} 
\begin{equation} \label{eq:omega}
\omega^t_\alpha(y) \,=\, \mathbb{P}_{\theta^{k-1}}(q_t = \alpha|y) \hspace{2.3cm}
\end{equation}
\end{subequations}
where $\mathbb{P}_{\theta^{k-1}}$ means the probability is computed under the current estimate $\theta^{k-1\,}$.

Thus, the E step (\ref{eq:estep}) amounts to implementing formulae (\ref{eq:nu}) and (\ref{eq:omega}). This is efficiently realised using Levinson's forward-Backward algorithm, recalled in Remark 2.

On the other hand, the M step (\ref{eq:mstep}) amounts to maximising the function (\ref{eq:qnew}), which is of the form (\ref{eq:lc}), and can be maximised as explained in Remark 3. \\[0.1cm]
$\diamond$ \textit{\textbf{Remark 2}\,:} in order to implement formulae (\ref{eq:estepprac}), it is helpful to use the forward variables $\Phi_t(\alpha)$ and backward variables $\mathrm{B}_t(\alpha)$, which satisfy
\begin{subequations} \label{eq:fb1}
\begin{equation} \label{eq:epractice1}
  \nu_{\alpha\beta}(y) \,=\, \sum^{T-1}_{t=1} \Phi_t(\alpha)\left(P_{\alpha\beta}f(y_{t+1}|\bar{y}_{\beta\,},\sigma_{\beta})\right)\mathrm{B}_{t+1}(\beta) \\[0.2cm]
\end{equation}
\begin{equation} \label{eq:epractice2}
  \omega^t_\alpha(y)\,=\, \Phi_t(\alpha)\,\mathrm{B}_t(\alpha) \phantom{fuckyouericmoulinesasss\,\,} \\[0.2cm] 
\end{equation}
\end{subequations}
and can be computed using Levinson's forward-backward algorithm, in terms of the forward recursion
\begin{subequations} \label{eq:fb2}
\begin{equation} \label{eq:forward}
  \Phi_{t+1}(\beta) \,=\,\Phi^{-1}_t\left[ \sum_{\;\alpha \in S}\Phi_t(\alpha)P_{\alpha\beta} \right]f(y_{t+1}|\bar{y}_{\beta\,},\sigma_{\beta}) \\[0.2cm]
\end{equation}
and of the backward recursion,
\begin{equation} \label{eq:backward}
  \mathrm{B}_t(\alpha) \,=\, \Phi^{-1}_{t+1}\,\sum_{\beta \in S}P_{\alpha\beta\,}\mathrm{B}_{t+1}(\beta)\,f(y_{t+1}|\bar{y}_{\beta\,},\sigma_{\beta}) \phantom{a\,} \\[0.2cm]
\end{equation}
where $\Phi_t$ is a normalising factor,
\begin{equation} \label{eq:normalisingfact}
\Phi_t \,=\, 
\sum_{\beta \in S}\left[ \sum_{\;\alpha \in S}\Phi_t(\alpha)P_{\alpha\beta} \right]f(y_{t+1}|\bar{y}_{\beta\,},\sigma_{\beta}) \\[0.1cm]
\end{equation}
and where (\ref{eq:forward}) is taken with the initial condition
\begin{equation} \label{eq:initialforward}
\Phi_{\scriptscriptstyle 1}(\alpha) \,=\, \frac{\mathstrut\pi_{\scriptscriptstyle 1}(\alpha)f(y_{\scriptscriptstyle 1}|\bar{y}_{\alpha\,},\sigma_\alpha)}{\mathstrut\sum_{\alpha\in S}
\mathstrut\pi_{\scriptscriptstyle 1}(\alpha)f(y_{\scriptscriptstyle 1}|\bar{y}_{\alpha\,},\sigma_\alpha)} \\[0.2cm]
\end{equation}
\end{subequations}
and (\ref{eq:backward}) is taken with the terminal condition $\mathrm{B}_{\scriptscriptstyle T}(\beta) = 1$.

The recursions (\ref{eq:forward}) and (\ref{eq:backward}) are efficient alternatives to Baum's forward-backward algorithm, introduced in\!\cite{rabiner}. 

They have the same computational complexity (roughly, the order of $T|S|^2$ where $|S|$ is the number of elements in $S$). On the other hand, they do not suffer from the same numerical instability, which makes Baum's algorithm rather difficult to use\!\cite{revisited}. \hfill$\blacksquare$

The proof of (\ref{eq:fb1}) and of (\ref{eq:fb2}) is discussed in Appendix \ref{app:fb}.

$\diamond$ \textit{\textbf{Remark 3}\,:} consider the task of maximising a function of the form (\ref{eq:lc}),
\begin{equation} \label{eq:Q}
\tiny{Q(\theta) = \sum_{\alpha,\beta \in S}\nu_{\alpha\beta}\log\left(P_{\alpha\beta}\right) \,+\, \sum_{\alpha \in S}\sum^T_{t=1} \omega^t_\alpha\log\,f(y_t|\bar{y}_{\alpha\,},\sigma_{\alpha})} \\[0.15cm]
\end{equation}
with $\nu_{\alpha\beta\,},\,\omega^t_\alpha > 0$ and $f(y|\bar{y},\sigma)$ given by (\ref{eq:f}), where the natural parameter $\eta(\sigma)$ is negative   -- if $\eta(\sigma)$ is positive, as in the von Mises-Fisher model, it is possible to absorb a minus sign into the function $D$. 

Then, any global maximiser of (\ref{eq:Q}), say $\hat{\theta} = (\hat{P}_{\alpha\beta\,},\hat{y}_\alpha\,,\hat{\sigma}_\alpha)$,\hfill\linebreak must satisfy 
$$
\hat{P}_{\alpha\beta} \,=\,\hat{P}_{\alpha\beta}(\nu)\;;\;
\hat{y}_\alpha \,=\, \hat{y}_{\alpha}(\omega) \;;\;
\eta(\hat{\sigma}_\alpha) = \hat{\eta}_\alpha(\omega) \\[0.15cm]
$$
where the following definitions are used \vspace{0.1cm}
\begin{subequations} \label{eq:estimz}
\begin{equation} \label{eq:estimp}
 \hat{P}_{\alpha\beta}(\nu) \,=\, \frac{\nu_{\alpha\beta}}{\mathstrut\sum_{\beta \in S}\nu_{\alpha\beta}} \\[0.1cm]
\end{equation}
\begin{equation} \label{eq:estimy}
 \hat{y}_{\alpha}(\omega) \,=\, \mathrm{argmin}_{y\in M}\,\sum^T_{t=1}\, \omega^t_\alpha\,D(y_t\,,\,y) \\[0.1cm]
\end{equation}
\begin{equation} \label{eq:estimsig}
 \hat{\eta}_\alpha(\omega) \,=\, \left(\psi^\prime\right)^{-1}\left(\frac{\sum^T_{t=1}\, \omega^t_\alpha\,D(y_t\,,\,\hat{y}_\alpha)}{\mathstrut \mathstrut\sum^T_{t=1}\,\omega^t_\alpha} \right) \\[0.1cm]
\end{equation}
\end{subequations}
Here, $\psi^\prime$ is the derivative of $\psi$ and $\left(\psi^\prime\right)^{-1}$ its reciprocal function. In particular, $\left(\psi^\prime\right)^{-1}$ is well-defined since $\psi$ is strictly convex. \hfill$\blacksquare$

The proof of Formulae (\ref{eq:estimz}) is discussed in Appendix \ref{app:estimz}. 

For the von Mises-Fisher model, the solution of the minimisation problem (\ref{eq:estimy}) is given by\!\cite{mardia}

\begin{subequations}
\begin{equation}\label{eq:estimyvmf}
  \hat{y}_\alpha(\omega) \,=\, U\left(\frac{\mathstrut\sum^T_{t=1}\,\omega^t_\alpha\,y_t}{\mathstrut \mathstrut\sum^T_{t=1}\,\omega^t_\alpha}\right) \\[0.1cm]
\end{equation}
with $U:\mathbb{R}^d - \lbrace 0 \rbrace \rightarrow \mathbb{R}^d$ the mapping $U(y) = y/\Vert y \Vert_{\scriptscriptstyle \mathbb{R}^d}$ where $\Vert\cdot\Vert_{\scriptscriptstyle \mathbb{R}^d}$ denotes the Euclidean norm.

For the Riemannian Gaussian model, (\ref{eq:estimy}) reads, after replacing (\ref{eq:rg2}), 
\begin{equation} \label{eq:estimyrg}
  \hat{y}_\alpha(\omega) \,=\,
  \mathrm{argmin}_{y\in M}\,\sum^T_{t=1}\, \omega^t_\alpha\,d^2(y_{t\,},y) \\[0.1cm]
\end{equation}
\end{subequations}
This means that $\hat{y}_\alpha(\omega)$ is the weighted Riemannian centre of mass of the observations $y_t\,$, and can be computed using any of the standard methods for computing Riemannian centres of mass\!\cite{said1}\!\!\cite{said2}.

Formula (\ref{eq:estimsig}) is easily evaluated. The function $\psi^\prime(\eta)$ being known, $\hat{\eta}_\alpha(\omega)$ is the unique solution of the equation
\begin{equation} \label{eq:equation}
 \psi^\prime(\hat{\eta}_\alpha(\omega)) \,=\, \frac{\sum^T_{t=1}\, \omega^t_\alpha\,D(y_t\,,\,\hat{y}_\alpha)}{\mathstrut \sum^T_{t=1}\, \omega^t_\alpha}
\end{equation}
Which can be found, to any accuracy, in exponential time, using the Newton-Raphson method, since the function $\psi(\eta)$ is strictly convex. 

It is also possible to find  $\hat{\eta}_\alpha(\omega)$ directly, by performing a search in a table which contains the values of the function $\psi^\prime(\eta)$. If the table is sufficiently large, this approach is successful within constant time, but has limited accuracy.








\vfill
\pagebreak

\section{Computer experiment} \label{sec:simus}
Here, a numerical experiment is considered, illustrating the hidden Markov chain model of Section \ref{sec:model}, and testing the validity of the EM algorithm of Section \ref{sec:algorithm}.

In this experiment, the hidden Markov chain $(q_t)$ takes values in a three-element set $S = \lbrace 1, 2, 3\rbrace$, and has initial distribution $\pi_{\scriptscriptstyle 0}(\alpha) = \delta_{\scriptscriptstyle 1}(\alpha)$, where  $\delta_{\scriptscriptstyle 1}(\alpha) \,=\,1$ if $\alpha = 1$ and $=0$ otherwise, and transition matrix $(P_{\alpha\beta})$,
\vspace{0.1cm}
\begin{equation} \label{eq:transitionalg}
(P_{\alpha\beta}) \,=\, \left(\begin{array}{ccc} 0.4 & 0.3 & 0.3 \\ 0.2 & 0.6 & 0.2 \\ 0.1 & 0.1 & 0.8 \end{array}\right)
\end{equation}
The outputs $(y_t)$ were generated from a Riemannian Gaussian model in the Poincar\'e disk. In other words, each $(y_t)$ has its values in $M = \lbrace z \in \mathbb{C}\,: |z| < 1\rbrace$, and the conditional density (\ref{eq:locationscale}) takes on the form (\ref{eq:rg3}), where\!\cite{vinberg}
\vspace{0.1cm}
\begin{subequations} \label{eq:poicnare1}
 \begin{equation} \label{eq:poincaredistance}
    d(y,z) \,=\, \mathrm{acosh}\!\left[ 1+ \frac{2|y\!-\!z|^2}{(1\!-\!|y|^2)(1\!-\!|z|^2)}\right] \\[0.2cm]
 \end{equation}
\begin{equation} \label{eq:poincarez}
  Z(\sigma) \,=\, (2\pi)^{\scriptscriptstyle 3/2}\,\sigma\,e^{\scriptscriptstyle \frac{\sigma^2}{\mathstrut 2}}\,\mathrm{erf}\left(\frac{\sigma}{\sqrt{2}}\right) \hspace{1.1cm}
\end{equation}
$\mathrm{erf}$ being the error function\!\cite{said2}.
\end{subequations}

The Poincar\'e disk is considered here, rather than a space $\mathcal{P}_{\scriptscriptstyle d}$ of positive-definite matrices, as in Section \ref{sec:model}, since it is much easier to visualise (indeed, it is just the inside of a unit disk in the plane). On the other hand, the Poincar\'e disk, with Riemannian distance (\ref{eq:poincaredistance}) is isometric to the space of $2\times 2$ symmetric positive-definite matrices, with unit determinant\!\cite{vinberg}.

Some background on the Riemannian geometry of the Poincar\'e disk is recalled in Remark 4.

The following sets of location parameters $\bar{y}_\alpha$ and scale parameters $\sigma_\alpha$ were used 
\begin{subequations} \label{eq:locscals}
\begin{equation} \label{eq:yalpha}
  \bar{y}_{\scriptscriptstyle 1} \,=\, 0 \hspace{0.11cm};\hspace{0.1cm}
  \bar{y}_{\scriptscriptstyle 2} \,=\, 0.82\,\mathrm{i} + 0.29  \hspace{0.11cm};\hspace{0.1cm}
  \bar{y}_{\scriptscriptstyle 3} \,=\, 0.82\,\mathrm{i} - 0.29
\end{equation}
\begin{equation} \label{eq:sigmalpha}
\sigma_{\scriptscriptstyle 1} \,=\, 0.1 \hspace{0.11cm};\hspace{0.1cm}
  \sigma_{\scriptscriptstyle 2} \,=\, 0.4 \hspace{0.11cm};\hspace{0.1cm}
  \sigma_{\scriptscriptstyle 3} \,=\, 0.4 
\end{equation}
\vspace{0.1cm}
where $\mathrm{i} = \sqrt{\mathstrut-1}$.
\end{subequations}

\begin{figure}
\begin{center}
\includegraphics[width=6.2cm]{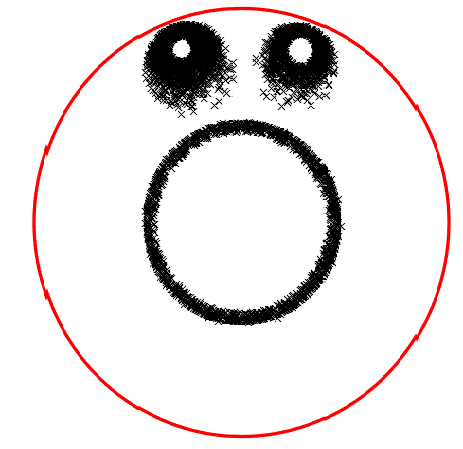}    
\caption{$10^4$ observations ; location-scale parameters (\ref{eq:locscals})} 
\label{fig:face1}
\end{center}
\end{figure}

With the hidden Markov model defined in this way, the observations $(y_t\,;t = 1,\ldots,T)$ where $T = 10000$ appeared as a ``face" pattern, within the Poincar\'e disk, shown in Figure \ref{fig:face1}. 

\begin{figure}
\begin{center}
\includegraphics[width=6.2cm]{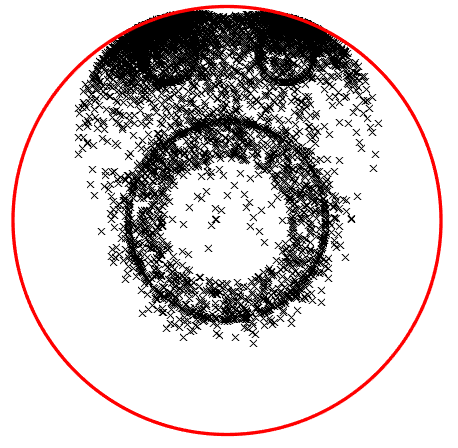}    
\caption{$10^4$ observations ; scale parameters $\times\; 5$} 
\label{fig:face2}
\end{center}
\end{figure}

If the scale parameters $\sigma_\alpha$ are all multiplied by $5$, the face becomes more blurred, since each conditional density (\ref{eq:locationscale}) has its dispersion multiplied by $5$. This can be seen in Figure \ref{fig:face2}.

In both figures \ref{fig:face1} and \ref{fig:face2}, the ``eyes" (centred at $\bar{y}_{\scriptscriptstyle 2}$ and $\bar{y}_{\scriptscriptstyle 3}$) appear smaller than the ``mouth", although are they have scale parameters $\sigma_{\scriptscriptstyle 2}$ and $\sigma_{\scriptscriptstyle 3}$ four times larger (see (\ref{eq:sigmalpha})\,!).  

This apparent contradiction is due to the fact that the Riemannian metric of the Poincar\'e disk is stretched near the boundary of the disk ---  see Formula (\ref{eq:metric}), Remark 4. Thus, while they appear smaller to us, objects near the boundary are, in fact, bigger when measured by this Riemannian metric.

Consider now the application of the EM algorithm of Section \ref{sec:algorithm}, to the estimation of the transition matrix (\ref{eq:transitionalg}) and location and scale parameters (\ref{eq:locscals}). 

\vfill
\pagebreak

$\diamond$ \textit{\textbf{Remark 4}\,:} the Poincar\'e disk is the set of complex numbers $z$ with $|z| < 1$, equipped with a conformal Riemannian metric, which turns it into a space of constant negative curvature\!\cite{vinberg}. 

This Riemannian metric measures the length of a plane vector $w$ (identified with a complex number), attached at the point $z$, using the Riemannian norm
\begin{equation} \label{eq:metric}
\Vert w \Vert^2_z \,=\, \left(\frac{4}{\kappa}\right)\frac{|w|^2}{\mathstrut (1-|z|^2)^2}
\end{equation}
where $\kappa > 0$ is a scaling factor, and $|w|^2 = w w^*$ is the squared modulus, or squared Euclidean norm (the star denotes the conjugate). 

It is called a conformal metric, since the Riemannian norm $\Vert w \Vert_z$ is proportional to the Euclidean norm $|w|$. The sectional curvature associated to this metric is constant, and equal to $-\kappa$.

In (\ref{eq:metric}), the conformal factor $(1-|z|^2)^{ -2}$ goes to infinity as $|z|$ approaches 1. Accordingly, the Riemannian norm $\Vert w \Vert_z$ of a vector $w$, with fixed Euclidean norm $|w|$, becomes arbitrarily large, when $z$ approaches the boundary of the disc. 

In other words, the Riemannian metric of the Poincar\'e disk is stretched near the boundary. \hfill$\blacksquare$ \\[0.2cm]
This algorithm was run $N_{\mathrm{mc}} = 20$ times, each time for $N_{\mathrm{em}} = 1000$ iterations. An individual run took about one hour, on a standard laptop.

All of these runs were carried out with the same initial guess $\theta^1$, but the observations $(y_t\,;t = 1,\ldots,T)$ were\hfill\linebreak generated anew for each run, leading to a different (in fact, random) value of the final estimate
\vspace{0.1cm}
$$
\theta^{N_{\mathrm{em}}} \,=\, (\hat{P}^{N_{\mathrm{em}}}_{\alpha\beta\,},\bar{y}^{N_{\mathrm{em}}}_\alpha\,,\sigma^{N_{\mathrm{em}}}_{\alpha})
$$
Then, the mean value of $(\hat{P}^{N_{\mathrm{em}}}_{\alpha\beta})$ was found to be
\vspace{0.1cm}
$$
\mathbb{E}_{\mathrm{mc}}(\hat{P}^{N_{\mathrm{em}}}_{\alpha\beta}) \,=\, 
\left(\begin{array}{ccc} 0.5 & 0.2 & 0.3 \\ 0.2 & 0.8 & 0 \\ 0.1 & 0 & 0.9 \end{array}\right)
$$
while the mean values of $\bar{y}^{N_{\mathrm{em}}}_\alpha$  and $\sigma^{N_{\mathrm{em}}}_\alpha$ turned out to be
\vspace{0.1cm}
$$
\begin{array}{llcl}
 \mathbb{E}_{\mathrm{mc}}(\bar{y}^{N_{\mathrm{em}}}_{\scriptscriptstyle 1}) \,=& - 0.02\,\mathrm{i} + 0.41 & ; &  \mathbb{E}_{\mathrm{mc}}(\sigma^{N_{\mathrm{em}}}_{\scriptscriptstyle 1}) \,=\, 0.5    \\[0.1cm]
  \mathbb{E}_{\mathrm{mc}}(\bar{y}^{N_{\mathrm{em}}}_{\scriptscriptstyle 2}) \,=&\phantom{-} 0.83\,\mathrm{i} + 0.30  & ; &  \mathbb{E}_{\mathrm{mc}}(\sigma^{N_{\mathrm{em}}}_{\scriptscriptstyle 2}) \,=\, 0.3  \\[0.1cm]
  \mathbb{E}_{\mathrm{mc}}(\bar{y}^{N_{\mathrm{em}}}_{\scriptscriptstyle 3}) \,=& \phantom{-}0.83\,\mathrm{i} - 0.29
& ; &  \mathbb{E}_{\mathrm{mc}}(\sigma^{N_{\mathrm{em}}}_{\scriptscriptstyle 3}) \,=\, 0.3  \\[0.2cm]
\end{array}
$$
Here, $\mathbb{E}_{\mathrm{mc}}$ is the ``Monte Carlo expectation", equal to the empirical average, taken over $N_{\mathrm{mc}}$ runs. 

The dispersion of $(\hat{P}^{N_{\mathrm{em}}}_{\alpha\beta})$, about its mean value, was measured by the maximum empirical variance of the entries $\hat{P}^{N_{\mathrm{em}}}_{\alpha\beta\,}$. This was found to be
$$
\max_{\alpha\beta}\,\mathbb{V}_{\mathrm{mc}}(\hat{P}^{N_{\mathrm{em}}}_{\alpha\beta}) \,=\, 0.06
$$
The dispersions of $\bar{y}^{N_{\mathrm{em}}}_\alpha$  and $\sigma^{N_{\mathrm{em}}}_\alpha$ were computed similarly
\vspace{0.1cm}
$$
\max_{\alpha}\,\mathbb{V}_{\mathrm{mc}}(\bar{y}^{N_{\mathrm{em}}}_\alpha) \,=\, 0.05 \;;\;\max_{\alpha}\,\mathbb{V}_{\mathrm{mc}}(\sigma^{N_{\mathrm{em}}}_\alpha) \,=\, 0.1
$$
where $\mathbb{V}_{\mathrm{mc}}$ is the ``Monte Carlo variance", equal to the empirical variance, taken over $N_{\mathrm{mc}}$ runs.

These numerical results show that the value of the final estimate $\theta^{N_{\mathrm{em}}}$ is rather stable, and is not affected very strongly by the randomness of the observations. Indeed, the dispersions of $\hat{P}^{N_{\mathrm{em}}}_{\alpha\beta}$ or of $\bar{y}^{N_{\mathrm{em}}}_\alpha$
and $\sigma^{N_{\mathrm{em}}}_\alpha$ appear small in comparison to their mean values.

Building on this remark, it is possible to infer that the discrepancy, between these mean values and the true values of the transition probabilities $P_{\alpha\beta\,}$,  and location and scale parameters $\bar{y}_\alpha$ and $\sigma_\alpha\,$, is due to the insufficient number of observations, $T = 10000$, rather than to a slow convergence of the EM algorithm. 

This is supported by the fact that each run of the EM algorithm was observed to converge to the final estimate $\theta^{N_{\mathrm{em}}}$ well before $N_{\mathrm{em}} = 1000$ iterations, producing nearly constant estimates $\theta^{k} \approx \theta^{N_{\mathrm{em}}}$ not later than $k = 800$. 

When the number of observations was increased to\hfill\linebreak $T = 50000$, an individual run, of $N_{\mathrm{em}} = 1000$ iterations of the algorithm, took much longer, over five hours.

For this reason, only one run of the algorithm was carried out. This was found to produce final estimates significantly closer to the true values of the transition matrix and location and scale parameters. 

However, for lack of time, a more detailed statistical analysis was not possible, in this case.

\section{Hidden Markov fields} \label{sec:fields}
Similar to hidden Markov chain models, hidden Markov field models involve a hidden stochastic process\footnote{a stochastic process is any family of jointly defined random variables\!\cite{kallenberg}.}$(q_z)$, which takes its values in a finite set $S$, observed through random outputs $(y_z)$, which belong to a homogeneous Riemannian manifold $M$. 

However, in the case of a hidden Markov field, the index $z$ is not a time variable (running through $t = 1,2,\ldots$), but a vertex of a finite unidrected graph $\mathcal{Z}$. Then, the stochastic process $(q_z\,;z \in \mathcal{Z})$ is called a Markov field if the following two properties are satisfied\!\cite{grimmett}
\begin{subequations} \label{eq:mrf}
\begin{equation} \label{eq:mrf1}
  \mathbb{P}(q_z = \alpha|q_w\,;w\neq z) \,=\, \mathbb{P}(q_z = \alpha|q_w\,;w\sim z) \\[0.1cm]
\end{equation}
where $w \sim z$ means $w$ and $z$ are adjacent, and
\begin{equation} \label{eq:mrf2}
 \mathbb{P}\left(\cap_{z \in \mathcal{Z}}\,\lbrace q_z = \varsigma(z)\rbrace\right)\,> 0 \\[0.1cm]
\end{equation}
for any function $\varsigma:\mathcal{Z} \rightarrow S$. Such a function is called a configuration of the field. 
\end{subequations}

Property (\ref{eq:mrf1}) says that the state $q_z$ at a vertex $z$ depends on states at the other vertices $q_w$ only through those $w$ which are adjacent to $z$. 

On the other hand, Property (\ref{eq:mrf2}) is required for the Hammersley-Clifford theorem\!\!\cite{grimmett}. This theorem implies that the probability distribution of a Markov field is a Gibbs distribution. 

For simplicity, it is here assumed that the graph $\mathcal{Z}$ is a square grid, where each vertex is adjacent to its immediate neighbors (this is a simple model of an image). In this case, the field $(q_z)$ is said to have a Gibbs distribution, if the probability mass function of the states $q_z$ takes on the form \vspace{0.1cm}
\begin{equation} \label{eq:gibbs}
p(q_z\,;z\in \mathcal{Z})\,\propto\,\prod_{z\in\mathcal{Z}}\exp\left[ V(q_z) + \frac{1}{2}\,\sum_{w\sim z} J(q_z\,,q_w)\right] \\[0.1cm]
\end{equation}
Here, $\propto$ indicates a missing normalising factor $W_{\scriptscriptstyle (V,J)}$. Computing $W_{\scriptscriptstyle (V,J)}$ requires summing over all possible\hfill\linebreak configurations of the field $q$ --- often, an impracticable task. 

The observed outputs $(y_z)$ depend on the hidden states $(q_z)$ exactly as in Section \ref{sec:model}. That is,
\begin{equation}\label{eq:locationscalebis}
  p(y_z|q_z = \alpha) \,=\, f(y_z|\bar{y}_{\alpha\,},\sigma_\alpha)
\\[0.1cm]
\end{equation}
where the function $f(y_z|\bar{y}_{\alpha\,},\sigma_\alpha)$ is of the form (\ref{eq:f}). Then, having access only to these outputs $(y_z)$, one hopes to recover the dependence structure of the Hidden field, given by the functions $V:S\rightarrow\mathbb{R}$ and $J:S\times S\rightarrow \mathbb{R}$, as well as the parameters $(\bar{y}_\alpha\,,\sigma_\alpha)$ which govern the $(y_z)$.

This problem generalises the hidden Markov chain model problem of Section \ref{sec:model}, which was addressed in Section \ref{sec:algorithm}. It is a problem of parameter estimation, for the unknown parameter $\theta = (V_\alpha\,,J_{\alpha\beta}\,,\bar{y}_\alpha\,,\sigma_\alpha)$, where $V_\alpha = V(\alpha)$ and $J_{\alpha\beta} = J(\alpha,\beta)$ for $\alpha,\beta \in S$. 

In the special case where the outputs $y_z$ belong to a Euclidean space $M = \mathbb{R}^d$, several EM algorithms were introduced in the image analysis literature, in order to address this problem\!\cite{meanf1}\!\!\!\cite{meanf2}. 

The extension of these algorithms, to observed outputs $y_z$ in any homogeneous Riemannian manifold $M$, is here discussed briefly.  

Each iteration of the EM algorithm includes an E step and an M step. Here, these take on the following form, where the notation is the same as in (\ref{eq:estep}) and (\ref{eq:mstep}) of Section \ref{sec:algorithm} --- unfortunately, due to lack of space, the following discussion does not include any proofs.\\[0.12cm]
$\bullet$ \textbf{E step}\,: compute the function \vspace{0.1cm}
\begin{equation} \label{eq:estepbis1}
Q(\theta|\theta^{k-1}) \,=\, Q(\bar{y}_\alpha\,,\sigma_\alpha) \,+\, Q(V,J) \\[0.1cm]
\end{equation}
in terms of the following formulae \vspace{0.1cm} 
\begin{subequations} \label{eq:estepbbis}
\begin{equation} \label{eq:estepbis2}
Q(\bar{y}_\alpha\,,\sigma_\alpha) \,=\, \sum_{\alpha \in S}\sum_{z\in\mathcal{Z}}\,\omega^z_\alpha(y)\log\,f(y_z|\bar{y}_{\alpha\,},\sigma_{\alpha})
\end{equation}
\begin{eqnarray} 
\nonumber Q(V,J) \,=\, -\log W(V,J) + \hspace{1.75cm}\\
\label{eq:estepbis3}  \sum_{\alpha \in S} \omega_\alpha(y)V_\alpha + \frac{1}{2}\,\sum_{\alpha,\beta \in S}\nu_{\alpha\beta}(y)J_{\alpha\beta}
\end{eqnarray}
\end{subequations}
which involve conditional probabilities and expectations,
\begin{subequations} \label{eq:omegabbis}
\begin{equation} \label{eq:omegabis1}
\omega^z_\alpha(y) \,=\, \mathbb{P}_{\theta^{k-1}}(q_z = \alpha|y) \hspace{0.1cm};\hspace{0.1cm}
\omega_\alpha(y) \,=\, \sum_{z \in \mathcal{Z}}\,\omega^z_\alpha(y)
\end{equation}
\begin{equation} \label{eq:omegabis2}
 \nu_{\alpha\beta}(y) \,=\, \sum_{w\sim z}\mathbb{P}_{\theta^{k-1}}(q_z = \alpha\,,q_{w}=\beta|y) \hspace{1.4cm}
\end{equation}  
\end{subequations}
$\bullet$ \textbf{M step}\,: maximise $Q(\theta|\theta^{k-1})$ to obtain new estimates,
$\theta^k = (\hat{V}_\alpha\,,\hat{J}_{\alpha\beta}\,,\hat{y}_\alpha\,,\hat{\sigma}_\alpha)$. These are given by,
\begin{subequations} \label{eq:estimbis}
\begin{equation} \label{eq:estimVJ}
(\hat{V},\hat{J}) \,=\, \mathrm{argmax}_{(V,J)}\, \langle \omega,V \rangle + \langle \nu,J\rangle - \Psi(V,J) \\[0.12cm]
\end{equation}
and, in the notation of Remark 3, \vspace{0.12cm}
\begin{equation} \label{eq:estimybis}
 \hat{y}_{\alpha} \,=\, \mathrm{argmin}_{\hat{y}\in M}\,\sum_{z\in \mathcal{Z}}\, \omega^z_\alpha(y)\,D(y_z\,,\,\hat{y}) \hspace{0.8cm}
\end{equation}
\vspace{0.12cm}
\begin{equation} \label{eq:estimsigbis}
 \eta(\hat{\sigma}_\alpha) \,=\, \left(\psi^\prime\right)^{-1}\left(\frac{\sum_{z\in \mathcal{Z}}\, \omega^z_\alpha(y)\,D(y_z\,,\,\hat{y}_\alpha)}{\mathstrut \sum_{z \in \mathcal{Z}}\,\mathstrut\omega^z_\alpha(y)} \right) \\[0.2cm]
\end{equation}
\end{subequations}
Here, $\Psi(V,J) = \log\,W(V,J)$, and 
$$
\langle\omega,V \rangle = \sum_{\alpha \in S} \omega_\alpha(y)V_\alpha \hspace{0.2cm};\hspace{0.2cm}
\langle \nu,J\rangle = \frac{1}{2}\sum_{\alpha,\beta \in S}\nu_{\alpha\beta}(y)J_{\alpha\beta}
$$ 
Moreover, the maximum in (\ref{eq:estimVJ}) is unique, since $\Psi(V,J)$ is a strictly convex function of $V = (V_\alpha)$ and $J= (J_{\alpha\beta})$, (when these are considered as $V \in \mathbb{R}^{|S|}$ and $J \in \mathbb{R}^{|S|\times|S|}$).  

$\diamond$ \textit{\textbf{Remark 5}\,:} in its form specified by (\ref{eq:estepbbis}) and (\ref{eq:estimbis}), the EM algorithm is not practically applicable. Indeed, both the normalising constant $W(V,J)$, and the conditional probabilities and expectations (\ref{eq:omegabbis}), are impossible to evaluate, outside of trivial situations.  

To circumvent this (fundamental) difficulty, various\hfill\linebreak approximate methods, for evaluating $W(V,J)$ and (\ref{eq:omegabbis}), have been proposed, such as the ones based on mean-field approximations\!\cite{meanf1}\!\!\cite{meanf2}. \hfill$\blacksquare$

It will be the aim of the authors' future work to adapt and apply these methods, to the general case where the observed outputs $y_z$ belong to a homogeneous Riemannian manifold.

\appendix

\section{Forward-backward variables} \label{app:fb}
Define the forward variables $\Phi_t(\alpha)$ and  backward variables $\mathrm{B}_t(\alpha)$ by\!\cite{revisited} --- for convenience, write $\mathbb{P}$ instead of $\mathbb{P}_{\theta^{k-1}}$
\vspace{0.1cm}
\begin{equation} \label{eq:fdef}
  \Phi_t(\alpha) \,=\, \mathbb{P}\left(q_t=\alpha|y_t\,,\ldots,y_{\scriptscriptstyle 1}\right) \phantom{abcd\,}\\[0.2cm]
\end{equation}
\begin{equation} \label{eq:bdef}
  \mathrm{B}_t(\alpha) \,=\, \frac{\mathbb{P}\left(y_{\scriptscriptstyle T}\,\ldots,y_{t+ 1}|q_t = \beta\right)}{\mathstrut p(y_{\scriptscriptstyle T}\,\ldots,y_{t+ 1}|y_t\,,\ldots,y_{\scriptscriptstyle 1})} \\[0.2cm]
\end{equation}
To see that these satisfy (\ref{eq:fb1}), the key is to make use of the fact that  $x_t =(y_{t\,},q_t)$ is a Markov chain, with transition probabilities\!\cite{hmmbook}\!\!\cite{moulines}
\begin{equation} \label{eq:xtransition}
  p(y_{t+1\,},q_{t+1} = \beta|y_{t\,},q_t = \alpha) \,=\, P_{\alpha\beta}\,f(y_{t+1}|\bar{y}_{\beta\,},\sigma_{\beta}) \\[0.2cm]
\end{equation}
For example, (\ref{eq:epractice1}) can be obtained by using Bayes rule to express each term in the sum (\ref{eq:nu}) 
\begin{equation} \label{eq:bayes}
\tiny{\mathbb{P}(q_t = \alpha\,,q_{t+1}=\beta|y) \,=\, \frac{p(y_{\scriptscriptstyle T}\,,\ldots,y_{\scriptscriptstyle 1}\,;q_t = \alpha\,,q_{t+1}=\beta)}{\mathstrut p(y_{\scriptscriptstyle T}\,,\ldots,y_{\scriptscriptstyle 1})} }\\[0.12cm]
\end{equation}
where the numerator is the joint probability density of $y_{\scriptscriptstyle T}\,,\ldots,y_{\scriptscriptstyle 1}$ on the event $(q_{t\,},q_{t+1}) =  (\alpha,\beta)$ --- this is a density with respect to the Riemannian volume measure of the product manifold $M^{\scriptscriptstyle T} = M \times \ldots \times M$.

This joint density can be written as a product,
\vspace{0.1cm}
$$
\begin{array}{rc}
p(y_{\scriptscriptstyle T}\,,\ldots,y_{t + 2}|\,y_{t+1\,},q_{t+1}=\beta\,,q_t = \alpha\,,y_{t}\,,\ldots,y_{\scriptscriptstyle 1}) &\times \\[0.12cm]
p(y_{t+1\,},q_{t+1} = \beta\,|\,q_t = \alpha\,,y_{t}\,,\ldots,y_{\scriptscriptstyle 1}) &\times \\[0.12cm]
\mathbb{P}(q_t = \alpha|\,y_{t}\,,\ldots,y_{\scriptscriptstyle 1}) &\times \\[0.12cm]
p(y_{t}\,,\ldots,y_{\scriptscriptstyle 1}) & \\
\end{array}
$$
by repeatedly applying (again!) Bayes rule. However, by the Markov property of $x_t =(y_{t\,},q_t)$, this product becomes
$$
\begin{array}{rc}
p(y_{\scriptscriptstyle T}\,,\ldots,y_{t + 2}|\,y_{t+1\,},q_{t+1}=\beta) &\times \\[0.12cm]
p(y_{t+1\,},q_{t+1} = \beta\,|\,q_t = \alpha) &\times \\[0.12cm]
\mathbb{P}(q_t = \alpha|\,y_{t}\,,\ldots,y_{\scriptscriptstyle 1}) &\times \\[0.12cm]
p(y_{t}\,,\ldots,y_{\scriptscriptstyle 1}) & \\
\end{array}
$$
but, using (\ref{eq:fdef}), (\ref{eq:bdef}) and (\ref{eq:xtransition}), this is the same as
$$
\begin{array}{rc}
\mathrm{B}_{t+1}(\beta) &\times \\[0.12cm]
P_{\alpha\beta}f(y_{t+1}|\bar{y}_{\beta\,},\sigma_{\beta}) &\times \\[0.12cm]
\Phi_t(\alpha) &\times \\[0.12cm]
p(y_{t}\,,\ldots,y_{\scriptscriptstyle 1}) & \\
\end{array}
$$
Therefore, replacing back into (\ref{eq:bayes}), it follows each term in the sum (\ref{eq:nu}) is given by
\vspace{0.1cm}
$$
\Phi_t(\alpha)\left(P_{\alpha\beta}f(y_{t+1}|\bar{y}_{\beta\,},\sigma_{\beta})\right)\mathrm{B}_{t+1}(\beta) \\[0.15cm]
$$
Thus, summing over $t = 1\,,\ldots, T-1$ yields (\ref{eq:epractice1}).

An analogous, and simpler, reasoning can be used to obtain (\ref{eq:epractice2}).

In addition to  (\ref{eq:fb1}), it should be proved that the forward and backward variables are given by (\ref{eq:fb2}).

To obtain (\ref{eq:forward}), note first that the initial condition (\ref{eq:initialforward}) follows directly from (\ref{eq:fdef}). In addition, also from (\ref{eq:fdef})
\vspace{0.1cm}
$$
\begin{array}{ll}
\Phi_{t+1}(\beta) &\propto\, \mathbb{P}\left(q_{t+1}=\beta|y_{t+1}\,,\ldots,y_{\scriptscriptstyle 1}\right) \\[0.12cm]
                     &\propto\, \mathbb{P}\left(y_{t+1\,},q_{t+1}=\beta|y_{t}\,,\ldots,y_{\scriptscriptstyle 1}\right) \\[0.12cm]
\end{array}
$$
where $\propto$ indicates a missing normalising factor. Thus, using the Markov property of $x_t =(y_{t\,},q_t)$, it is seen that, up to normalisation, $\Phi_{t+1}(\beta)$ is equal to 
$$
\begin{array}{ll}
& \sum_{\alpha \in S}p(y_{t+1\,};q_{t+1}=\beta|q_t = \alpha) \times \mathbb{P}(q_t=\alpha|y_{t}\,,\ldots,y_{\scriptscriptstyle 1}) \\[0.5cm]
=& \sum_{\alpha \in S}P_{\alpha\beta}\,f(y_{t+1}|\bar{y}_{\beta\,},\sigma_{\beta}) \times \Phi_t(\alpha)
\end{array}
$$
This immediately gives (\ref{eq:forward}), after noting that the missing normalising factor can be found by summing over $\beta \in S$, as in (\ref{eq:normalisingfact}). 

The proof of (\ref{eq:backward}) is similar to that of  (\ref{eq:forward}), and is here omitted to avoid repetition.

\section{Proof of Formulae (\ref{eq:estimz})} \label{app:estimz}
Recall the function $Q(\theta)$ of (\ref{eq:Q}),
\vspace{0.1cm}
$$
\tiny{Q(\theta) = \sum_{\alpha,\beta \in S}\nu_{\alpha\beta}\log\left(P_{\alpha\beta}\right) \,+\, \sum_{\alpha \in S}\sum^T_{t=1} \omega^t_\alpha\log\,f(y_t|\bar{y}_{\alpha\,},\sigma_{\alpha})} \\[0.15cm]
$$
$\diamond$ \textit{\textbf{Proof of (\ref{eq:estimp})}\,:} only the first term in the expression of $Q(\theta)$ depends on $(P_{\alpha\beta})$. Therefore, this term can be maximised separately. 

Denote this first term by $Q(P)$. Then, maximising $Q(P)$ is equivalent to maximising
$$
\tilde{Q}(P) \,=\,\sum_{\alpha,\beta \in S} \hat{P}_{\alpha\beta}(\nu)\log\left(P_{\alpha\beta}\right)
$$ 
where $\hat{P}_{\alpha\beta}(\nu)$ is given by (\ref{eq:estimp}).

Let $\hat{P}_{\alpha\beta} = \hat{P}_{\alpha\beta}(\nu)$. For each $\alpha \in S$, let $\hat{P}_\alpha = (\hat{P}_{\alpha\beta}\,;\beta \in S)$ and $P_\alpha = (P_{\alpha\beta}\,;\beta \in S)$. Both $\hat{P}_\alpha$ and $P_\alpha$ are probability distributions on $S$. If $\mathrm{KL}(\hat{P}_\alpha\Vert P_\alpha)$ is the Kullback-Leibler divergence between $\hat{P}_\alpha$ and $P_\alpha\,$, then\!\cite{cover}
$$
\tilde{Q}(P) \,=\, \sum_{\alpha \in S} \left\lbrace H(\hat{P}_\alpha) \,-\, \mathrm{KL}(\hat{P}_\alpha\Vert P_\alpha) \right\rbrace
$$
where $H(\hat{P}_\alpha) = \sum_{\beta \in S}\, \hat{P}_{\alpha\beta}\log(\hat{P}_{\alpha\beta})$ is the entropy of $\hat{P}_\alpha$ (which does not depend on $(P_{\alpha\beta})$). 

Therefore, $\tilde{Q}(\theta)$ reaches its maximum when each each $\mathrm{KL}(\hat{P}_\alpha\Vert P_\alpha)$ reaches its minimum. But this happens when $\mathrm{KL}(\hat{P}_\alpha\Vert P_\alpha) = 0$ at $P_{\alpha\beta} =\hat{P}_{\alpha\beta\,}$, as required. \hfill$\blacksquare$ \\[0.2cm] 
$\diamond$ \textit{\textbf{Proof of (\ref{eq:estimy}) and (\ref{eq:estimsig})}\,:} denote the second term in the expression of $Q(\theta)$ by $Q(\bar{y},\sigma)$ and consider the task of maximising $Q(\bar{y},\sigma)$.

This can be carried out by maximising each one of the functions
$$
Q(\bar{y}_\alpha\,,\sigma_\alpha) \,=\, \sum^T_{t=1} \omega^t_\alpha\log\,f(y_t|\bar{y}_{\alpha\,},\sigma_{\alpha}) \\[0.13cm]
$$
with respect to $\bar{y}_\alpha$ and $\sigma_\alpha\,$. From (\ref{eq:f}),
$$
Q(\bar{y}_\alpha\,,\sigma_\alpha) \,=\, \sum^T_{t=1} \omega^t_\alpha\left\lbrace\, \eta(\sigma_\alpha)D(y_t\,,\bar{y}_\alpha) - \psi(\eta(\sigma_\alpha))\right\rbrace \\[0.13cm]
$$
By first maximising over $\bar{y}_\alpha$ and recalling that $\eta(\sigma_\alpha)$ is negative, it is seen than any maximiser $\hat{y}_\alpha$ must be equal to $\hat{y}_\alpha(\omega)$ given by (\ref{eq:estimy}). 

Then, it remains to maximise (over $\sigma_\alpha$) the function
$$
Q(\sigma_\alpha)\,=\,\sum^T_{t=1} \omega^t_\alpha\left\lbrace\, \eta(\sigma_\alpha)D(y_t\,,\hat{y}_\alpha) - \psi(\eta(\sigma_\alpha))\right\rbrace
$$
This is a strictly concave function of $\eta(\sigma_\alpha)$, since $\psi(\eta)$ is a strictly convex function. Thus, it is possible to maximise by differentiating with respect to $\eta(\sigma_\alpha)$ and setting the derivative to zero. This yields the equation 
$$
\sum^T_{t=1} \omega^t_\alpha\left\lbrace\, D(y_t\,,\hat{y}_\alpha) - \psi^\prime(\eta(\hat{\sigma}_\alpha))\right\rbrace \,=\, 0
$$
for the maximiser $\hat{\sigma}_\alpha\,$. This equation is equivalent to (after dividing both sides by $\sum^T_{t=1}\,\omega^t_\alpha$)
\vspace{0.1cm}
$$
\psi^\prime(\eta(\hat{\sigma}_\alpha)) \,=\, \frac{\sum^T_{t=1}\, \omega^t_\alpha\,D(y_t\,,\,\hat{y}_\alpha)}{\mathstrut \sum^T_{t=1}\, \omega^t_\alpha} \\[0.1cm]
$$
which is the same as equation (\ref{eq:equation}), whose unique solution is $\eta(\hat{\sigma}_\alpha) = \hat{\eta}_\alpha(\omega)$ given by (\ref{eq:estimsig}). \hfill$\blacksquare$ 



\end{document}